\documentclass[10pt,leqno,twoside]{amsart}

\usepackage{amsmath,amssymb,amsfonts,amsthm}
\usepackage{fancyhdr}
\usepackage{amssymb,amsbsy}
\usepackage[left=4cm,right=4cm,top=4cm,bottom=4cm,foot=1.7cm]{geometry}

\newtheorem{theorem}{\indent Theorem}[section]

\newtheorem{lemma}{\indent Lemma}[section]

\newtheorem{remark}{\indent Remark}[section]

\newcommand\Var{\operatorname{\textrm{Var}}}

\title{Weights and degrees in a random graph model based on 3-interactions} 

\author{\'Agnes Backhausz, Tam\'as F. M\'ori}

\dedicatory{\upshape 
Department of Probability Theory and Statistics,\\ 
E\"otv\"os Lor\'and University\\ 
P\'azm\'any P.~s. 1/C, H-1117 Budapest, Hungary\\                         
\textit{E-mail address:} \texttt{agnes@cs.elte.hu,
 moritamas@ludens.elte.hu}
} 

\thanks{The European Union and the European Social Fund have
provided financial support to the project under the grant agreement no.\ 
T\'AMOP 4.2.1./B-09/KMR-2010-0003.}

\date{4 June 2012}

\begin{document}

\begin{abstract}
In a random graph model introduced in \cite{BM} we 
give the joint asymptotic distribution of weights and degrees and
prove scale-free property for the model. Moreover, we determine 
the asymptotics of the maximal weight and the maximal degree.
\end{abstract}

\maketitle

\section{Introduction}

Many random graph models have been invented recently for modelling large 
networks like the internet or social networks \cite{BA, durrett}. 
Considering degree distributions, real life networks look quite different
from classical (i.e., Erd\H os--R\'enyi type) random graphs. Motivated
by this observation, in several models \cite{BA, cf} the evolution of
the graph is driven by the actual degrees. However, in real-world
networks larger groups and cliques may also interact and this  
has a relevant effect on the evolution. Therefore models based on cliques   
or groups of vertices may be of particular interest \cite{[Mo09], nastos}.  

Motivated by that, in \cite{BM} we introduced a random graph model 
with dynamics based on interactions of three vertices. In our model
vertices taking part in an interaction together have larger chance to
participate together again, thus it is a kind of preferential
attachment structure, while in the models referred to earlier there 
is no possibility to keep track of the number of steps where the members 
of a given group get new edges together.

In our model vertices, edges and triangles have nonnegative weights,
increasing randomly in discrete time steps. The weight is the number
of interactions that the vertex, pair of vertices or triplets
have participated in. In \cite{BM} we investigated the limit of the
ratio of vertices of a given weight and proved that they almost surely
exist and decay polynomially. This is the scale-free property of the
model. We also determined the asymptotics of the weight of a given
vertex.  

This time we will deal with degrees. This is the number of vertices 
having interacted with a given one. This is different from the 
weight of the vertex, which is the total number of interactions. 
We determine the joint asymptotic distribution of weights and degrees; 
prove scale-free property for degrees; finally, we give 
the asymptotics of the degree of a given vertex, the maximal weight,
and the maximal degree.

\section{The model}

We start with a single triangle. This has initial weight $1$, and all
its three edges have weight $1$. Vertices, edges, and triangles will
have nonnegative integer-valued weights, which increase according
to the random evolution of the graph.  

At each step three vertices will interact. There are two
possibilities. With probability $p$, independently of the past, a new
vertex is added, which then interacts with two already existing
vertices. Otherwise three old vertices interact. We will need $0<p\leq
1$.  

Assume that in the $n$th step a new vertex is added to the graph.  
With probability $r$, independently of the past, the
choice is done according to the ``preferential attachment'' rule, that is, 
an edge is chosen with probability proportional to its weight, then its 
endpoints are selected. With probability $1-r$ two distinct old
vertices are chosen uniformly at random. 

Then the new vertex interacts with the two selected vertices. 
This means that the triangle they form comes to existence with initial
weight $1$, and we increase the weights of all three edges of the
$3$-interaction by $1$. This is the end of the step where a new vertex is
generated.    

With probability $1-p$ three of the old vertices will interact. In
such a step we have two choices again. With probability $q$ each
triangle is selected with probability proportional to its weight.
Otherwise, with probability $1-q$, three distinct vertices will 
be chosen at random, uniformly, i.e., each triplet with the same
probability. This choice is also independent of the past.  

In each case, having selected the three vertices to interact,
we draw the edges of the triangle that are not present yet. Then
the weight of the triangle is increased by $1$, as well as the weights
of the three sides of the triangle. 

Now we define the weights of vertices. The weight of a vertex is the
sum of the weights of the triangles that contain it. Note that this is
just the half of the sum of weights of edges from it, because whenever
a vertex takes part in an interaction, the first sum is increased
by $1$, and the latter one is increased by $2$.   

Our model is parametrized by the triplet of probabilities $(p,q,r)$.

This construction was introduced in \cite{BM}, where the following
properties were proved.  

The ratio of vertices of weight $w$ converges to $x_w$ almost surely
as $n\rightarrow \infty$, where
\begin{equation}\label{wrec}
x_1=\frac{1}{\alpha+\beta+1}\,,\quad
x_w=\frac{\alpha(w-1)+\beta}{\alpha w+\beta+1}\,x_{w-1},
\end{equation}
hence we have
\begin{equation}\label{wasympt}
x_w\sim\frac{\Gamma\bigl(1+\tfrac{\beta+1}{\alpha}\big)}
{\alpha\,\Gamma\bigl(1+\tfrac{\beta}{\alpha}\bigr)}\,
w^{-\bigl(1+\tfrac{1}{\alpha}\bigr)},\quad\text{as }w\to\infty
\end{equation}
\cite[Theorem 3.1]{BM}.

We have also studied the rate of growth of the weight of a fixed
vertex.

It is clear that the weights of the vertices of the starting triangle
are interchangeable, therefore it is not necessary to deal with all
the three. Let them be labelled by $-2$, $-1$, and $0$. The further
vertices get labels $1$, $2$, etc, in the order they are added to the
graph. Let $D[n,j]$ and $W[n,j]$ denote the degree and the weight of vertex
$j$ after step $n$, provided it exists. Otherwise let these quantities
be equal to zero. Obviously, vertex $j\ge 1$ cannot exist before step $j$.

According to Theorem 4.1 of \cite{BM}, for $j\ge 0$ fixed we have 
\begin{equation}\label{W[n,j]}
W[n,j]\sim \zeta_jn^{\alpha}\quad\text{almost surely as }n\to\infty,
\end{equation}
where $\zeta_j$ is a positive random variable.

In the sequel we will denote by $\mathcal F_n$ the $\sigma$-field
generated by the first $n$ steps, and by $V_n$ the number of vertices
after the $n$th step. Thus $V_0=3$. Furthermore, let $\mathbb
I(\,\cdot\,)$ be defined as $1$ if the condition within the brackets
holds, otherwise let it be $0$.  

\section{Asymptotic joint distribution of degree and weight}

We denote the number of vertices of weight $w$ and degree $d$ after
$n$ steps by $X[n,d,w]$. When a vertex is born, its initial weight is
one, and its initial degree is two. When it takes part in an
interaction, its weight is increased by one, while its degree may not
change (if it is already connected to the other two interacting
vertices), or may increase by one or two. Thus  $X[n,d,w]>0$ can occur
only for pairs of integers $d,w$ with $1\leq w$ and $2\leq d
\leq 2w$. 

The following theorem is about the almost sure convergence
of the ratio of vertices of weight $w$ and degree $d$.
\begin{theorem} \label{hat1.1}
Given integers $1\leq w$ and $2\leq d\leq 2w$ we have 
\[
\frac{X[n,d,w]}{V_n}\rightarrow x_{d,w}
\]
almost surely as $n\rightarrow \infty$, where the limits $x_{d,w}$ 
are positive numbers satisfying the following recurrence equation. 
\begin{equation*}
\begin{split}
x_{2,1}&=\frac{1}{\alpha+\beta+1}\,,\\
x_{d,w}&=\frac{1}{\alpha w+\beta+1}\big[\alpha_1 (w-1)
x_{d,w-1}+\alpha_2 (w-1) x_{d-1,w-1}+\beta x_{d-2,w-1} \big]
\end{split}
\end{equation*}
for $w\geq 2$, where 
\begin{equation}\label{harmaxe1}
\alpha_1=(1-p)q, \quad \alpha_2=\frac{2pr}{3},\quad 
\alpha=\alpha_1+\alpha_2,\quad  
\beta=\frac{1}{p}[2(1-r)+3(1-p)(1-q)].
\end{equation}
\end{theorem}

\proof We compute the conditional expectation of $X[n,d,w]$ with
respect to the $\sigma$-algebra $\mathcal F_{n-1}$. Note that if an
old vertex interacts with a new one, its degree must increase. On
the other hand, if we choose vertices with probabilities proportional
to certain weights, then no new edges are born between old vertices. 

Having built the graph in $n$ steps we consider a fixed vertex with
degree $d$ and weight $w$. For simplicity we denote by $V=V_n$ the
number of vertices after $n$ steps. Then $(d,w)$ can increase 
\begin{itemize}
\item by $(1,1)$ with probability $pr\,\dfrac{2w}{3(n+1)}$ (new
  vertex, preferential attachment);
\item by $(1,1)$ with probability $p(1-r)\,\dfrac{d}{\binom V2}\,$,
and\\
by $(2,1)$ with probability $p(1-r)\,\dfrac{V-d-1}{\binom V2}$ (new
vertex, uniform selection); 
\item by $(0,1)$ with probability $(1-p)q\,\dfrac{w}{n+1}$ (old
  vertices, preferential attachment);
\item by $(0,1)$ with probability $(1-p)(1-q)\,\dfrac{\binom d2}
{\binom V3}\,$,\\
by $(1,1)$ with probability $(1-p)(1-q)\,\dfrac{d(V-d-1)}{\binom
  V3}\,$,  and\\
by $(2,1)$ with probability
$(1-p)(1-q)\,\dfrac{\binom{V-d-1}{2}}{\binom V3}$ (old vertices,
uniform selection). 
\end{itemize}

Now it is easy to see that the probability that a vertex of weight $w$
takes part in the interaction at step $n$ is the following (see
also \cite{BM}). 
\begin{equation}\label{hae1.0}
p\Bigg[r\,\frac{2w}{3n}+(1-r)\,\frac{2}{V_{n-1}}\Bigg]+
(1-p)\Bigg[q\,\frac wn+(1-q)\,\frac{3}{V_{n-1}}\Bigg]= 
\frac{\alpha w}{n}+\frac{\beta p}{V_{n-1}}\,;
\end{equation} 
this is independent of the degree of the vertex.

For $d=2$ and $w=1$ we also have to take the new vertex into account:
a new vertex is born with probability $p$, and its degree is surely
$2$, while its weight is surely $1$. 

Summing up, we obtain the conditional expectation of $X[n,d,w]$ in the
following form.
\begin{equation}\label{xdwe1}
\begin{split}
\mathbb E(X[n,d,w]&\mid\mathcal F_{n-1})= 
X[n-1,d,w]\left(1-\frac{\alpha w}{n}-\frac{\beta p}{V_{n-1}}\right)\\
&+X[n-1,d,w-1](1-p)\left[q\frac{w-1}{n}+(1-q)\frac{\binom{d}{2}}
{\binom{V_{n-1}}{3}}\right]\\
&+X[n-1,d-1,w-1]\,p\left[r\frac{2(w-1)}{3n}+(1-r)\frac{d}{\binom{V_{n-1}}{2}}
\right]\\
&+X[n-1,d-1,w-1](1-p)(1-q)\frac{d(V_{n-1}-1-d)}{2\binom{V_{n-1}}{3}}\\
&+X[n-1,d-2,w-1]\times\\
&\quad\times\left[p(1-r)\frac{V_{n-1}-d-1}{\binom{V_{n-1}}{2}}+
(1-p)(1-q)\frac{\binom{V_{n-1}-1-d}{2}}{\binom{V_{n-1}}{3}}\right]\\
&+p\,\mathbb I(d=2,w=1).
\end{split}
\end{equation}

Introduce the normalizing sequence
\[
c[n,w]=\prod_{i=1}^{n-1} \left(1-\frac{\alpha w}{i}-
\frac{\beta p}{V_{i-1}}\right)^{\!-1},\quad n\geq 1,\, w\geq 1.
\]
At each step a new vertex is born with probability $p$ independently
of the past. Hence the law of large numbers can be applied to the
number of vertices, yielding that 
\begin{equation}\label{hae1}
V_n=pn+o\left(n^{1/2+\varepsilon}\right)
\end{equation} 
a.s., for all $0<\varepsilon<\frac 12$. This implies that
\begin{multline*}
\log c[n,w]=\sum_{i=1}^{n-1} -\log \Biggl(1-\frac{\alpha w}{i}-
\frac{\beta}{i+o\left(i^{1/2+\varepsilon}\right)}\Biggr)\\
=\sum_{i=1}^{n-1}\left(\frac{\alpha w}{i}+\frac{\beta}{i}
+o\bigl(i^{-3/2+\varepsilon}\bigr)\right)=
(\alpha w+\beta)\sum_{i=1}^{n-1}\frac{1}{i}+O(1)
\end{multline*} 
a.s., where the error term is convergent as
$n\rightarrow\infty$. Therefore  
\begin{equation}\label{hae2}
c[n,w]\sim a_w n^{\alpha w+\beta}
\end{equation}
a.s. as $n\rightarrow\infty$, where $a_w$ is a positive random
variable.  
 
For $n\geq 1,\, w\geq 1, d\geq 2$ set $Z[n,d,w]=c[n,w]X[n,d,w]$. 
Multiplying both sides of \eqref{xdwe1} by $c[n,w]$ one can see that
$\bigl(Z[n,d,w],\,\mathcal F_n\bigr)$ is a nonnegative submartingale
for all fixed integers $w, d$. 
  
Consider the Doob decomposition $Z[n,d,w]=M[n,d,w]+A[n,d,w]$, where
$M[n,d,w]$ is a zero mean martingale, and $A[n,d,w]$ is a predictable
increasing process. 
\begin{align*}
M[n,d,w]&=\sum_{i=1}^n\Bigl(Z[i,d,w]-\mathbb E\bigl(Z[i,d,w]\bigm|
\mathcal F_{i-1}\bigr)\Bigr),\\
A[n,d,w]&=\mathbb EZ[1,d,w]+\sum_{i=2}^n\Bigl(\mathbb E\bigl(
Z[i,d,w]\bigm|\mathcal F_{i-1}\bigr)-Z[i-1,d,w]\Bigr).
\end{align*}

Let us give an upper bound on the conditional variance of the
martingale part. Recall
that $c[i,w]$ is $\mathcal F_{i-1}$-measurable, and since there is
only one interaction at each step, the increment of $X$ can not be
greater than three. Using \eqref{hae2} we get that  
\begin{equation}\label{hafok1}
\begin{split}
B[n,d,w]&=\sum_{i=1}^n \Var(Z[i,d,w]\,|\,\mathcal F_{i-1})=
\sum_{i=1}^n c[i,w]^2\Var(X[i,d,w]\,|\,\mathcal F_{i-1})\\
&=\sum_{i=1}^n c[i,w]^2\Var\bigl(X[i,d,w]-X[i-1,d,w]\bigm|\mathcal
F_{i-1}\bigr)\\ 
&\leq \sum_{i=1}^n c[i,w]^2\, \mathbb E \left(\left.\bigl(X[i,d,w]-
X[i-1,d,w]\bigr)^2\right\vert \mathcal F_{i-1}\right)\\
&\leq 9 \sum_{i=1}^n c\left[i, w\right]^2=
O\left( n^{2(\alpha w+\beta)+1}\right).
\end{split}
\end{equation}
Let us apply Proposition VII-2-4 of Neveu \cite{[Ne75]} with
$f(t)=\sqrt t \log t$. We obtain that  
\begin{equation}\label{mart}
M[n,d,w]=o\bigl(B[n,d,w]^{1/2}\log B[n,d,w]\bigr)=o\bigl(n^{\alpha
 w+\beta+1}\bigr).
\end{equation}
We will later use this estimation to show that $Z[n,d,w]\sim A[n,d,w]$.

Note that $X[n,d,w]=0$ if $2\leq d\leq 2w$ does not hold. Hence
$x_{d,w}=0$ in all these cases.

We apply induction on $w$. If the weight of a vertex is equal to $1$,
then it could not participate in any interactions except the first
one, when it was born. Thus its degree must be equal to two. Therefore
$X[n,d,1]$ is zero for $d\neq 2$, and it is  the number of vertices of
weight $1$ for $d=2$. In the case $w=1$ the proposition follows  
from \eqref{wrec}.

Suppose that the statement holds for all weights less than $w$, and
for all possible degrees $2\leq d\leq 2w$. Let us compute the
asymptotics of $A[n,d,w]$. We start from \eqref{xdwe1}.

\begin{align*}
A[n,d,w]=\mathbb E Z&[1,d,w]\\
&+\sum_{i=2}^{n}c[i,w]X[i-1,d,w-1](1-p)\left[q\frac{w-1}{i}+
(1-q)\frac{\binom{d}{2}}{\binom{V_{i-1}}{3}}\right]\\
&+c[i,w]X[i-1,d-1,w-1]\,p\left[r\frac{2(w-1)}{3i}+(1-r)
\frac{d}{\binom{V_{i-1}}{2}}\right]\\
&+c[i,w]X[i-1,d-1,w-1](1-p)(1-q)\frac{d(V_{i-1}-1-d)}
{2\binom{V_{i-1}}{3}}\\
&+c[i,w]X[i-1,d-2,w-1]\times\\
&\quad\times \left[p(1-r)\frac{V_{i-1}-d-1}{\binom{V_{i-1}}{2}}
+(1-p)(1-q)\frac{\binom{V_{i-1}-1-d}{2}}{\binom{V_{i-1}}{3}}\right]\\
&+c[i,w]p\,\mathbb I(d=2,w=1).
\end{align*}
  
Using the induction hypothesis, the asymptotics of $V_n$ in
\eqref{hae1} and the regular variation of the normalizing constants
$c[n,w]$ with exponent $\alpha w+\beta$ in equation \eqref{hae2}, we
can compute the asymptotics of $A[n,d,w]$, leaving out all terms that
are of smaller order of magnitude than others.
\begin{align*}
A[n,d,w]&\sim\sum_{i=2}^{n}c[i,w]\,pi\,x_{d,w-1}\,(1-p)q\frac{w-1}{i}\\
&\qquad+c[i,w]\,pi\,x_{d-1,w-1}\,pr\frac{2(w-1)}{3i}\\
&\qquad+c[i,w]\,pi\,x_{d-2,w-1}\left[2(1-r)
+\frac{3(1-p)(1-q)}{p}\right]\\
&\sim p\sum_{i=2}^{n} a_w i^{\alpha w+\beta}\Bigg[(1-p)q(w-1)x_{d,w-1}
+pr\frac{2(w-1)}{3}\,x_{d-1,w-1}\\ 
&\qquad+\left(2(1-r)+\frac{3(1-p)(1-q)}{p}\right)x_{d-2,w-1}\Bigg]\\
&\sim p\,\frac{a_w n^{\alpha w+\beta+1}}{\alpha w+\beta+1}
\Bigg[(1-p)q(w-1)x_{d,w-1}+pr\frac{2(w-1)}{3}x_{d-1,w-1}\\
&\qquad+\left(2(1-r)+\frac{3(1-p)(1-q)}{p}\right) x_{d-2,w-1}\Bigg]. 
\end{align*} 
If $d$ and $w$ satisfy $2\leq d\leq 2w$, then there is at
least one term on the right-hand side that is positive due to the
induction hypothesis. Thus $M[n,d,w]=o\bigl(A[n,d,w]\bigr)$, therefore
$Z[n,d,w]\sim A[n,d,w]$ holds almost surely as $n\to 
\infty$. Dividing by the normalizing constants $c[n,w]$ we get that   
$X[n,d,w]\sim x_{d,w} pn$, from which
\[
\frac{X[n,d,w]}{V_n}\rightarrow x_{d,w}
\]
a.s. as $n\rightarrow \infty$, where
\begin{align*}
x_{d,w}&=\frac{1}{\alpha w+\beta+1}\biggl[(1-p)q(w-1)x_{d,w-1}+
pr\,\frac{2(w-1)}{3}\,x_{d-1,w-1}\\
&\hspace{35mm}+\left(2(1-r)+\frac{3(1-p)(1-q)}{p}\right) x_{d-2,w-1}\biggr]\\
&=\frac{1}{\alpha w+\beta+1}\,\left[\alpha_1 (w-1) x_{d,w-1}+\alpha_2 (w-1)
x_{d-1,w-1}+\beta x_{d-2,w-1}\right], 
\end{align*}
with
\[
\alpha_1=(1-p)q,\quad\alpha_2=\frac{2pr}{3}\,,\quad 
\beta= 2(1-r)+\frac{3(1-p)(1-q)}{p}\,.
\]

By this the induction step is completed. Moreover, as we noted before,
$x_{d,w}>0$ holds for $2\leq d\leq 2w$. \qed  

\begin{remark} 
The explicit solution of the recurrence equation in
Theorem \ref{hat1.1} can be given in the following form. 

For $w\ge 1$ set 
\[
c_w=(\alpha w+\beta+1)(\alpha(w-1)+\beta+1)\dots(\alpha+\beta+1).
\]
Let $S_n(0)=1$, and for $1\le k\le n$ define
\[
S_n(k)=\sum_{1\le i_1<i_2<\dots<i_k\le n}i_1 i_2\dots i_k.
\]
Then
\begin{equation}\label{explicit}
x_{d,w}=\frac{1}{c_w}\sum_{k=1}^w S_{w-1}(w-k)\binom{w-k}{d-2k} 
\alpha_1^{w-d+k}\alpha_2^{d-2k}\beta^{k-1},
\quad 1\leq w,\ 2\leq d \leq 2w.
\end{equation}
In other words, $x_{d,w}$ is equal to the coefficient of $z^{d-2}$ in
the expression
\[
\frac{1}{c_w}\prod_{i=1}^{w-1}\bigl(i(\alpha_1+\alpha_2z)+\beta z^2\bigr).
\]
This is not hard to derive, and even easier to check; however, it does
not seem to be very convenient for determining the asymptotics of $x_{d,w}$
as $d$ or $w$ tends to infinity. We rather choose another method for it
in the next section.  
\end{remark}

\section{Construction of the two dimensional limit distribution}
Let $W$ be a positive integer valued random variable with
distribution $\mathbb P(W=w)=x_w$, $w=1,2,\dots$\,. In addition, let
$\xi_1\equiv 2$, and let the random variables $\xi_2,\xi_3,\dots$ be
independent of each other and of $W$ too; moreover, 
\begin{gather*}
\mathbb P(\xi_w=0)=\frac{\alpha_1(w-1)}{\alpha(w-1)+\beta}\,,\quad
\mathbb P(\xi_w=1)=\frac{\alpha_2(w-1)}{\alpha(w-1)+\beta}\,,\\ 
\mathbb P(\xi_w=2)=\frac{\beta}{\alpha(w-1)+\beta}\,. 
\end{gather*}
Define the partial sums $S_w=\xi_1+\dots+\xi_w$.
\begin{theorem}\label{const}
\[
\mathbb P(S_W=d,\,W=w)=x_{d,w},\quad 1\le w,\ 2\le d\le 2w.
\]
\end{theorem}
\proof
Clearly, $\mathbb P(S_W=2,\,W=1)=\mathbb P(W=1)=x_1=x_{2,1}$. In
addition we have
\begin{multline*}
\mathbb P(S_W=d,\,W=w)=\mathbb P(S_w=d,\,W=w)\\
\begin{aligned}
&=\mathbb P(S_w=d)
\mathbb P(W=w)\\
&=\Bigl[\mathbb P(S_{w-1}=d)\mathbb P(\xi_w=0)+
\mathbb P(S_{w-1}=d-1)\mathbb P(\xi_w=1)+\\
&\qquad+\mathbb P(S_{w-1}=d-2)\mathbb P(\xi_w=2)\Bigr]\Bigl[
\mathbb P(W=w-1)\,\frac{\alpha(w-1)+\beta}{\alpha w+\beta+1}\Bigr]\\
&=\mathbb P(S_{w-1}=d,\,W=w-1)\,\frac{\alpha_1(w-1)}{\alpha
  w+\beta+1}+\\ 
&\qquad+\mathbb P(S_{w-1}=d-1,\,W=w-1)\,\frac{\alpha_2(w-1)}{\alpha
  w+\beta+1}+\\
&\qquad+\mathbb P(S_{w-1}=d-2,\,W=w-1)\,\frac{\beta}{\alpha
  w+\beta+1}\,,
\end{aligned}
\end{multline*}
thus the probabilities $\mathbb P(S_W=d,\,W=w)$ and the limits
$(x_{d,w})$ satisfy the same recursion. 
\qed

This, combined with Theorem \ref{hat1.1}, implies that the empirical joint
distribution of degree and weight after step $n$ converges almost
surely in total variation norm to the distribution of $(S_W,\,W)$.
As a corollary we obtain that the asymptotic weight distribution
$(x_w)$ is just the marginal of the joint distribution $(x_{d,w})$, 
that is, $x_w=x_{1,w}+\dots+x_{2w,w}$.

\begin{theorem}\label{bivarlim}
Suppose both $\alpha_1$ and $\alpha_2$ are positive. Then
\[
x_{d,w}=x_w\cdot\frac{\alpha}{\sqrt{2\pi\alpha_1\alpha_2 w}}
\biggl(\exp\Bigl(-\frac{(\alpha d-\alpha_2 w)^2}
{2\alpha_1\alpha_2w}\Bigr)+O\bigl(w^{-1/2}\bigr)\biggr),
\]
as $d$ and hence $w$ tend to infinity; and the term $O$ in the
remainder is uniform in $d$.
\end{theorem}
\proof
\[
\mathbb E\xi_w=\frac{\alpha_2(w-1)+2\beta}{\alpha(w-1)+\beta}=
\frac{\alpha_2}{\alpha}+\frac{(\alpha+\alpha_1)\beta}
{\alpha(\alpha(w-1)+\beta)}\,,
\]
hence $\mathbb ES_w=\dfrac{\alpha_2}{\alpha}\,w+O(\log w)$. Similarly,
\[
\Var(\xi_w)=\frac{\alpha_1\alpha_2}{\alpha^2}+O\Bigl(\frac{1}{w}\Bigr)\,,
\quad \Var(S_w)=\frac{\alpha_1\alpha_2}{\alpha^2}\,w+O(\log w),
\]
as $w\to\infty$. The proof can be completed by applying the local
limit theorem (Theorem VII.1.5 in \cite{Petrov}) to $S_w$. Its
conditions are satisfied, namely,
\[
\liminf_{w\to\infty}\frac{1}{w}\Var(S_w)>0,\quad
\limsup_{w\to\infty}\frac{1}{w}\sum_{j=1}^w\left|\xi_j-\mathbb E\xi_j
\right|^3<\infty.
\]
Hence we have
\begin{equation}\label{llt}
\sup_{d\in\mathbb Z}\left|\sqrt{\Var(S_w)}\,\mathbb P(S_w=d)-
\frac{1}{\sqrt{2\pi}}\,\exp\left(-\frac{(d-\mathbb ES_w)^2}
{2\Var(S_w)}\right)\right|=O\left(\frac{1}{\sqrt{w}}\right).
\end{equation}
It easily follows that in \eqref{llt} $\mathbb ES_w$ can be replaced
with a term differing from it by $O\biggl(\sqrt{\dfrac{w}{\log w}}\,
\biggr)$, and $\Var(S_w)$ with a term differing by
$O\bigl(\sqrt{w}\bigr)$.   
\qed
\bigskip

From Theorem \ref{const} one can derive the asymptotics of the other
marginal distribution $u_d=\sum_{w\ge d/2}x_{d,w}$. Clearly, $u_d$ is
the a.s.\ limit of the proportion of vertices with degree $d$.  
\begin{theorem}
\[
u_d\sim\frac{\Gamma\bigl(1+\tfrac{\beta+1}{\alpha}\big)}
{\alpha_2\,\Gamma\bigl(1+\tfrac{\beta}{\alpha}\bigr)}\,
\Bigl(\frac{\alpha}{\alpha_2}\,d\Bigr)^{-\bigl(1+\tfrac{1}{\alpha}\bigr)},
\quad\text{as }d\to\infty.
\]
\end{theorem}
\proof
Let
\begin{gather*}
f=\frac{\alpha}{\alpha_2}\,d,\quad H=\{w: f-f^{1/2+\varepsilon}\le 
w\le f+f^{1/2+\varepsilon}\},\\
H^-=\{w: w<f-f^{1/2+\varepsilon}\},\quad 
H^+=\{w: w>f+f^{1/2+\varepsilon}\},
\end{gather*}
with some $0<\varepsilon<1/6$.

By Hoeffding's well-known exponential inequality (Theorem 2 of
\cite{H}) for $w\in H^-$ we have
\begin{multline*}
\mathbb P(S_w\ge d)\le \mathbb P\left(S_w-\mathbb ES_w\ge 
d-\frac{\alpha_2}{\alpha}\,w-O(\log w)\right)\\
\le\exp\Biggl(-\frac{\bigl(d-\tfrac{\alpha_2}{\alpha}w-O(\log w)
\bigr)^2}{2w}\Biggr)=
\exp\Biggl(-\Bigl(\frac{\alpha_2}{\alpha}\Bigr)^{\!2}
\frac{\bigl(f-w-O(\log w)\bigr)^2}{2w}\Biggr).
\end{multline*}
Here in the numerator $\bigl(f-w-O(\log w)\bigr)^2\ge 
f^{1+2\varepsilon}-O\bigl(f^{1/2+\varepsilon}\log f\bigr)$, and in the 
denominator $w\le f$. Hence 
\[
\mathbb P(S_w\ge d)\le\exp\biggl(-\frac{\alpha_2^2f^{2\varepsilon}}
{2\alpha^2}+o(1)\biggr),
\]
thus we have
\begin{equation}\label{also}
\mathbb P(S_W=d,\,W\in H^-)\le \bigl(1+o(1)\bigr)f\,
\exp\biggl(-\frac{\alpha_2^2 f^{2\varepsilon}}{2\alpha^2}\biggr)=
o\Bigl(f^{-\bigl(1+\tfrac{1}{\alpha}\bigr)}\Bigr).
\end{equation}

The case of $w\in H^+$ can be treated similarly. 
\begin{multline*}
\mathbb P(S_w\le d)\le \mathbb P\left(S_w-\mathbb ES_w\le 
d-\frac{\alpha_2}{\alpha}\,w\right)\\
\le\exp\biggl(-\frac{\bigl(\tfrac{\alpha_2}{\alpha}w-d\bigr)^2}{2w}
\biggr)\le\exp\biggl(-\Bigl(\frac{\alpha_2}{\alpha}\Bigr)^{\!2}
\frac{(w-f)^2}{2w}\biggr).
\end{multline*}
This time we use the estimate 
\[
2(w-f)\ge f^{1/2+\varepsilon}+w-f\ge f^{1/2+\varepsilon}+
(w-f)^{1/2+\varepsilon}\ge w^{1/2+\varepsilon}
\]
in the numerator, obtaining
\[
\mathbb P(S_w\le d)\le\exp\biggl(-\frac{\alpha_2^2w^{2\varepsilon}}
{8\alpha^2}\biggr).
\]
Hence
\begin{equation}\label{felso}
\mathbb P(S_W=d,\,W\in H^+)\le \sum_{w>f}
\exp\biggl(-\frac{\alpha_2^2 w^{2\varepsilon}}{8\alpha^2}\biggr)=
o\Bigl(f^{-\bigl(1+\tfrac{1}{\alpha}\bigr)}\Bigr).
\end{equation}

Finally, for $w\in H$
\begin{multline*}
\frac{(\alpha d-\alpha_2w)^2}{2\alpha_1\alpha_2w}=\frac{\alpha_2(f-w)^2}
{2\alpha_1w}=\frac{\alpha_2(f-w)^2}{2\alpha_1f}\Bigl(1+O\bigl(f^{-1/2+
\varepsilon}\bigr)\Bigr)\\
=\frac{\alpha_2(f-w)^2}{2\alpha_1f}+O\bigl(f^{-1/2+3\varepsilon}\bigr),
\end{multline*}
consequently
\[
x_{d,w}\sim\frac{\Gamma\bigl(1+\tfrac{\beta+1}{\alpha}\big)}
{\alpha\,\Gamma\bigl(1+\tfrac{\beta}{\alpha}\bigr)}\,
f^{-\bigl(1+\tfrac{1}{\alpha}\bigr)}\cdot
\frac{\alpha}{\sqrt{2\pi\alpha_1\alpha_2 f}}
\exp\biggl(-\frac{\alpha_2(w-f)^2}{2\alpha_1f}\biggr),
\]
as $d\to\infty$ and $w\in H$. Since
\[
\sum_{w\in H} \frac{\alpha}{\sqrt{2\pi\alpha_1\alpha_2 f}}
\exp\biggl(-\frac{\alpha_2(w-f)^2}{2\alpha_1f}\biggr)\to
\int\limits_{-\infty}^{+\infty}\frac{\alpha}
{\sqrt{2\pi\alpha_1\alpha_2}}\exp\biggl(-\frac{\alpha_2t^2}
{2\alpha_1}\biggr)dt=\frac{\alpha}{\alpha_2},
\]
we obtain that
\begin{equation}\label{kozep}
\mathbb P(S_W=d,\,W\in H)\sim\frac{\Gamma\bigl(1+\tfrac{\beta+1}
{\alpha}\big)}{\alpha_2\,\Gamma\bigl(1+\tfrac{\beta}{\alpha}\bigr)}\,
f^{-\bigl(1+\tfrac{1}{\alpha}\bigr)}.
\end{equation}

The proof is completed by \eqref{also}, \eqref{felso}, and
\eqref{kozep} combined. 
\qed

\section{Maximal weight, maximal degree}

In this section our goal is to determine the asymptotics of the
maximum of the weights as the number of steps tends to infinity. 

Let $I[n,j]$ denote the indicator of the event $\{W[n,j]\ge
1\}$. Moreover, we denote by $J[n,j]$ the indicator of the event that
vertex $j$ is born at step $n$, that is, $J[n,j]=I[n, j]-I[n-1,j]$.  

We fix $j$, and examine the process $W[n,j]$ as $n$
increases. 

In the first lemma we find martingales that we will use later in the
proofs. Then we prove that the maximal weight grows at the
same pace as the weight of any fixed vertex does, see \eqref{W[n,j]}.

Let $j,\,k,\,\ell$ be fixed integers, $0\le j\le\ell$, $1\le k$, and
let us introduce the sequences
\[
b[n,k]=\prod_{i=1}^{n}\Bigl(1+\frac{\alpha k}{i}\Bigr)^{\!-1},\quad 
d[n,k,j]=\sum_{i=1}^{n-1} b[i+1,k]\frac{\beta p}{V_i}
\binom{W[i,j]+k-1}{k-1}\,, 
\]
with $\alpha$, $\beta$ defined in \eqref{harmaxe1}. Note that $b[n,k]$
is deterministic, while $d[n,k,j]$ is random, but $\mathcal
F_{n-1}$-measurable for all $k$ and $j$. Moreover, we have 
\begin{equation}\label{harmaxe2}
b[n,k]\sim b_k n^{-k\alpha},
\end{equation}
with $b_k>0$, as $n\rightarrow \infty$.  

\begin{lemma}\label{hal1}
Let
\[
Z[n,k,j]=b[n,k]\binom{W[n,j]+k-1}{k}-d[n,k,j],
\]
then $\bigl(Z[n,k,j]I[\ell,j],\,\mathcal F_n\bigr)$, $n\ge\ell$, is a
martingale.
\end{lemma}
\proof 
Assuming that vertex $j$ already exists after step $\ell$, 
the probability that it participates in an interaction at step
$n+1\ge\ell$ is equal to $\frac{\alpha W\left[n,j\right]}
{n+1}+\frac{\beta p}{V_n}$. This implies that, for arbitrary positive 
integers $k,\,\ell,\,n$ we have 
\begin{multline*}
\mathbb E\left(\left. \binom{W[n+1,j]+k-1}{k}I[\ell,j]\,\right\vert
\mathcal F_n\right)=I[\ell,j]\binom{W[n,j]+k-1}{k}\\
+I[\ell,j]\left[\frac{\alpha W[n,j]}{n}+\frac{\beta p}{V_n}\right]
\left[\binom{W[n,j]+k}{k}-\binom{W[n,j]+k-1}{k}\right]\\
=I[\ell,j]\binom{W[n,j]+k-1}{k}+I[\ell,j]\left[\frac{\alpha W[n,j]}{n}
+\frac{\beta p}{V_n}\right]\binom{W[n,j]+k-1}{k-1}\\
=I[\ell,j]\binom{W[n,j]+k-1}{k}\left( 1+\frac{\alpha k}{n}\right)
+I\left[\ell,j\right]\frac{\beta p}{V_n}\binom{W\left[n,j\right]+k-1}{k-1}.
\end{multline*}
Multiplying both sides by $b[n+1,k]$ we get by definition that 
\begin{multline*}
\mathbb E\left(\left. b[n+1,k]\binom{W[n+1,j]+k-1}{k}I[\ell,j]\right
\vert \mathcal F_n\right)\\
=I[\ell,j]\left[b[n,k]\binom{W[n,j]+k-1}{k}+\frac{\beta p}{V_n}
\binom{W[n,j]+k-1}{k-1} b[n+1,k]\right]\\
=I[\ell,j]\left[b[n,k]\binom{W[n,j]+k-1}{k} -d[n,k,j]+d[n+1,k,j]\right],
\end{multline*}
which completes the proof of the lemma, since $d[n+1,k,j]$ is $\mathcal
F_n$-measurable. \qed

\begin{lemma}\label{harmaxl2}
For arbitrary fixed integers $k\ge 0$ and $1\le m\leq n$ define
\[
S[m,n,k]=\sum_{j=m}^n \mathbb E\left[b[n,k]\binom{W[n,j]+k-1}{k}  
I[n,j]\right].
\]
Then
\[
S[m,n,k]\le C_k\sum\limits_{j=m}^n j^{-k\alpha}
\]
with a positive constant $C_k$ only depending on $k$.
\end{lemma}
\proof

We prove this by induction on $k$. 

For $k=0$ obviously
\begin{equation*}
S[m,n,0]=b[n,0]\sum_{j=m}^n \mathbb E\bigl(I[n,j]\bigr)=
\sum_{j=m}^n \mathbb P\bigl(W[n,j]\ge 1\bigr)\le n-m+1.
\end{equation*}

Suppose that the statement of the lemma holds for $k-1$. By Lemma
\ref{hal1} we know that $Z[n,k,j]I[\ell,j]$ is a martingale,  
hence its expectation does not depend on $n$. The difference of
martingales is also a martingale, thus we have the same with
$J[\ell,j]$. Decomposing $I[n,j]$ into the sum of terms $J[\ell,j]$ we 
obtain that   
\begin{multline*}
S[m,n,k]=
\sum_{j=m}^n \mathbb E \left(\sum_{\ell=j}^n b[n,k]
\binom{W[n,j]+k-1}{k} J[\ell,j]\right)\\
=\sum_{j=m}^n\mathbb E\left(\sum_{\ell=j}^n\bigl(Z[n,k,j]+
d[n,k,j]\bigr)J[\ell,j]\right)\\
=\sum_{j=m}^n\mathbb E\left(\sum_{\ell=j}^n\bigl(Z[\ell,k,j]+
d[n,k,j]\bigr)J[\ell,j]\right)\\
=\mathbb E \left(\sum_{j=m}^n \sum_{l=j}^n 
\bigl(b[\ell,k]+d[n,k,j]-d[\ell,k,j]\bigr)J[\ell,j]\right).
\end{multline*}

Let us split $S[m,n,k]$ into two parts: $S[m,n,k]=S_1[m,n,k]+S_2[m,n,k]$,
where
\begin{multline}\label{S1}
S_1[m,n,k]=\mathbb E\left(\sum_{j=m}^n \sum_{\ell=j}^n b[\ell,k]
J[\ell,j]\right)\leq \sum_{j=m}^n b[j,k]\mathbb E\left(\sum_{\ell=j}^n
J[\ell,j]\right)\\
\le\sum_{j=m}^n b[j,k]=b_k \sum_{j=m}^n j^{-k\alpha}
\bigl(1+o(1)\bigr)\le C_k^\prime \sum_{j=m}^n j^{-k\alpha}.
\end{multline}

For the second part we have 
\begin{multline*}
S_2[m,n,k]=\mathbb E\left(\sum_{j=m}^n \sum_{\ell=j}^n 
\bigl(d[n,k,j]-d[\ell,k,j]\bigr)J[\ell,j]\right)\\
=\mathbb E\left( \sum_{j=m}^n \sum_{\ell=j}^n \sum_{i=\ell}^{n-1}b[i+1,k]
\,\frac{\beta p}{V_i}\binom{W[i,j]+k-1}{k-1}J[\ell,j]\right)\\ 
=\mathbb E\left(\sum_{i=m}^{n-1}b[i+1,k]\,\frac{\beta p}{V_i} 
\sum_{j=m}^i \sum_{\ell=j}^{i} \binom{W[i,j]+k-1}{k-1}J[\ell,j]
\right)\\
=\mathbb E\left( \sum_{i=m}^{n-1}\frac{b[i+1,k]}{b[i,k-1]}\,
\frac{\beta p}{V_i} \sum_{j=m}^i b[i,k-1]\binom{W[i,j]+k-2}{k-1}
\;\frac {W[i,j]+k-1}{W[i,j]}\;I[i,j]\right).
\end{multline*}
Since $\dfrac{W[i,j]+k-1}{W[i,j]}\le k$, we get that
\[
S_2[m,n,k]\le k\sum_{i=m}^{n-1} \frac{b[i+1,k]}{b[i,k-1]}\,\mathbb E
\left( \frac{\beta p}{V_i} \sum_{j=m}^i b[i,k-1]
\binom{W[i,j]+k-2}{k-1} I[i,j]\right).
\]

For the expectation in the right-hand side we give upper bounds on the
events $\lbrace V_i<(p/2)i\rbrace$ and $\lbrace V_i\geq (p/2)i\rbrace$
separately. Remember that $\mathbb I(\,\cdot\,)$ denotes the indicator
of the event in brackets. From the induction hypothesis we obtain that  
\begin{multline*}
\mathbb E\left(\frac{\beta p}{V_i}\,\mathbb I\Bigl(V_i\ge\frac{pi}{2}\Bigr)
\sum_{j=m}^i b[i,k-1]\binom{W[i,j]+k-2}{k-1}I[i,j]\right)\\
\le \frac{2 \beta}{i}S[m,i,k-1]\le 2\beta C_{k-1}\,\frac{1}{i}
\sum_{j=m}^i j^{-(k-1)\alpha}. 
\end{multline*}

By the Hoeffding bound $\mathbb P\bigl(V_i<(p/2)i\bigr)\le 
e^{-\varepsilon i}$ with $\varepsilon >0$ only depending on
$p$. Using that $V_i\geq 3$ and the trivial bound on the weights  
we get that 
\begin{multline*}
\mathbb E\left(\frac{\beta p}{V_i}\,\mathbb I\Bigl(V_i<\frac{pi}2\Bigr)
\sum_{j=m}^i  b[i,k-1]\binom{W[i,j]+k-2}{k-1} I[i,j]\right) \\
\le\frac{\beta p}{3}\,\mathbb P\Bigl(V_i< \frac{pi}2\Bigr)
\sum_{j=m}^i  b[i,k-1]\binom{i+k-2}{k-1}\\
=O \left(e^{-\varepsilon i} i^{-(k-1)\alpha} i^{k-1} i \right)=
o\left(\frac{1}{i} \sum_{j=m}^i j^{-(k-1)\alpha}\right)
\end{multline*}
uniformly in $m$. Finally,
\[
\frac{b[i+1,k]}{b[i,k-1]}=O\bigl(i^{-\alpha}\bigr).
\]

Putting all these together we obtain that 
\begin{multline}\label{S2}
S_2[m,n,k]\leq C_k'' \sum_{i=m}^n i^{-1-\alpha}
\sum_{j=m}^i j^{-(k-1)\alpha}
=C_k'' \sum_{j=m}^n j^{-(k-1)\alpha}\sum_{i=j}^n i^{-1-\alpha}\\
\le C_k'''\sum_{j=m}^n j^{-k\alpha}.
\end{multline}

We can complete the proof by combining \eqref{S1} and \eqref{S2}.\qed

\bigskip

Next we characterize the growth rate of the maximal weight in the
graph. Let $\mathcal W_n=\max\{W[n,j]:-2\le j\le n\}$, the maximal
weight after $n$ steps.

\begin{theorem}\label{maxw}
$\mathcal W_n\sim\mu\,n^{\alpha}$ almost surely as $n\to\infty$,
where $\mu$ is a finite and positive random variable, namely,
$\mu=\sup\{\zeta_j:j\ge -2\}$, with $\zeta_j$ defined in \eqref{W[n,j]}.  
\end{theorem}

\proof 
For $1\le m\le n$ define $M[m,n]=\max\{W[n,j]: -2\le j<m\}$. By 
\eqref{W[n,j]} it is obvious that
\[
\lim_{n\to\infty}n^{-\alpha}M[m,n]=\max\{\zeta_j: -2\le j<m\}
\]
with probability $1$. All we have to do is to show that
\begin{equation}\label{cel}
\lim_{m\to\infty}\limsup_{n\to\infty}n^{-\alpha}
\bigl(\mathcal W_n-M[m,n]\bigr)=0.
\end{equation}

From the proof of Lemma \ref{hal1} it follows that the process 
\[
b[n,k]\binom{W[n,j]+k-1}{k} I[\ell,j],\quad n\ge\ell,
\] 
is a submartingale, hence, the same holds for
\[
b[n,k]\binom{W[n,j]+k-1}{k}=b[n,k]\binom{W[n,j]+k-1}{k}I[n,j],
\quad n\ge j.
\]
Being the maximum of an increasing number of submartingales, the
process  
\[
b[n,k]\binom{\mathcal W_n-M[m,n]+k-1}{k},\quad n\ge m,
\]
is also a submartingale. In addition,  
\begin{equation}\label{becsles1}
\mathbb E \left[ b[n,k]\binom{\mathcal W_n-M[m,n]+k-1}{k}\right ]
\leq S[m,n,k]\le C_k\sum_{j=m}^n j^{-k\alpha}
\end{equation}
by Lemma \ref{harmaxl2}. Since
\begin{equation}\label{becsles2}
\Bigl([b[n,1]\bigl(\mathcal W_n-M[m,n]\bigr)\Bigr)^{\!k}\le
k!\,\frac{b[n,1]^k}{b[n,k]}\,b[n,k]\binom{\mathcal W_n-M[m,n]+k-1}{k},
\end{equation}
the nonnegative submartingale $b[n,1]\bigl(\mathcal W_n-M[m,n]\bigr)$
is bounded in $L_k$ whenever $k\alpha>1$. Thus it is convergent with
probability $1$, and also in $L_k$, for every $k\ge 1$. Moreover, by
\eqref{becsles1} and \eqref{becsles2}  we have
\[
\mathbb E\left(\lim_{n\to\infty}n^{-\alpha}
\bigl(\mathcal W_n-M[m,n]\bigr)\right)^{\!k}\le k!\,\frac{C_k}{b_k}
\sum_{j=m}^\infty j^{-k\alpha}.
\]
From this the monotone convergence theorem gives
\[
\mathbb E\left(\lim_{m\to\infty}\lim_{n\to\infty}n^{-\alpha}
\bigl(\mathcal W_n-M[m,n]\bigr)\right)^{\!k}=0
\]
if $k>1/\alpha$, proving \eqref{cel}. \qed

We finally present the asymptotics of the
maximal degree as the number of steps tends to infinity. 
First we will study the growth of the degree of a fixed vertex.

\begin{theorem}\label{D[n,j]}
For $j=0,1,\dots $ we have
\[
D[n,j]\sim\frac{\alpha_2}{\alpha}\,\zeta_j\,n^{\alpha}
\]
almost surely, as $n\to\infty$, where $\zeta_j$ is a positive random
variable, defined in \eqref{W[n,j]}.
\end{theorem}

\proof
Starting from the specification of the ways the degree and the weight
of a fixed vertex can grow, we can write
\begin{multline}\label{fokfelt}
\mathbb E\bigl(I[k,j]D[n+1,j]\bigm|\mathcal F_n\bigr)=
I[k,j]\Biggl(D[n,j]+pr\frac{2W[n,j]}{3(n+1)}\\
+p(1-r)\frac{2V_n-D[n,j]-2}
{\binom{V_n}{2}}+3(1-p)(1-q)\frac{V_n-D[n,j]-1}{\binom{V_n}{2}}\Biggr)\\
=I[k,j]\biggl(D[n,j]+\alpha_2 \frac{W[n,j]}{n+1}+R_n\biggr),
\end{multline}
if $k\le n$, where $0\le R_n\le\dfrac{2p\beta}{V_n}$. 

Introduce $\xi_n=I[k,j]\bigl(D[n,j]-D[n-1,j]\bigr)$, then
$0\le\xi_n\le 2$, hence by Corollary VII-\textbf{2}-6 of \cite{[Ne75]}
and equation \eqref{W[n,j]} we obtain 
\[
\sum_{i\le n}\xi_i\sim\sum_{i\le n}\mathbb E(\xi_i\mid\mathcal F_{i-1})
=\sum_{i\le n}\Bigl(\alpha_2\frac{W[i-1,j]}{i}+R_{i-1}\Bigr)\sim
\frac{\alpha_2}{\alpha}\,\zeta_j\,n^{\alpha},
\]
a.e. on the event $\{W[k,j]\ge 1\}$. 
Thus $D[n,j]\sim\frac{\alpha_2}{\alpha}\,\zeta_j\,n^{\alpha}$ on that
event. Since we know that $\lim_{k\to\infty}W[k,j]=\infty$, we have
\[
\mathbb P\biggl(\bigcup_k \{W[k,j]\ge 1\}\biggr)=1,
\]
completing the proof.\qed

\bigskip

From Theorems \ref{maxw} and \ref{D[n,j]} the asymptotic behaviour of
the maximal degree immediately follows. 
\begin{theorem}\label{maxd}
Let $\mathcal D_n$ denote the maximal degree in the graph after $n$ steps.
Then
\[
\mathcal D_n\sim \frac{\alpha_2}{\alpha}\mu\, n^{\alpha}
\]
almost surely as $n\to\infty$, where $\mu$ is the finite and positive
random variable defined in Theorem \ref{maxw}.
\end{theorem}
\proof 
By the trivial bound $D[n,j]\le 2W[n,j]$ we obtain
\begin{multline*}
\max\{D[n,j]:-2\le j<m\}\le\mathcal D_n\\
\le\max\{D[n,j]:-2\le j<m\}+\max\{2W[n,j]:m\le j\le n\}.
\end{multline*}
Multiplying by $n^{-\alpha}$ and letting $n\to\infty$ we get
\begin{multline*}
\frac{\alpha_2}{\alpha}\max\{\zeta_j:-2\le j<m\}\le\liminf_{n\to\infty}
n^{-\alpha}\mathcal D_n\le\limsup_{n\to\infty}n^{-\alpha}\mathcal D_n\\
\le\frac{\alpha_2}{\alpha}\max\{\zeta_j:-2\le j<m\}+
2\lim_{n\to\infty}n^{-\alpha}\bigl(\mathcal W_n-M[m,n]\bigr).
\end{multline*}
By \eqref{cel} both sides tend to $\mu\alpha_2/\alpha$ as 
$m\to\infty$.\qed

\vspace{1cm}

\noindent\textbf{\'Agnes Backhausz}\\
Department of Probability Theory and Statistics, 
 E\"otv\"os Lor\'and University, Budapest, Hungary\\
{\tt agnes@cs.elte.hu}\
\bigskip 

\noindent\textbf{Tam\'as F. M\'ori}\\
Department of Probability Theory and Statistics,  
 E\"otv\"os Lor\'and University, Budapest, Hungary\\
 {\tt moritamas@ludens.elte.hu}

\end{document}